\documentclass{article}
\usepackage{amsfonts}
\usepackage{amsmath}
\usepackage{amssymb}
\usepackage[english]{babel}

\DeclareMathOperator{\vol}{Vol}
\DeclareMathOperator{\diam}{diam}
\DeclareMathOperator{\inter}{int}

\begin {document}

\begin{center}
CONTINUALITY OF SET OF BILIPSCHITZ CLASSES 

IN EUCLIDEAN SPACE

\medskip

{\it A.~Magazinov}
\end{center}

{\bf Introduction}
\medskip

This paper is devoted to studying biLipschitz equivalence of Delone sets.

Let $M$ be a metric space with distance $d_M(x,y)$. Denote by $B_{\rho}(x)$ and $B_{\rho}^\circ(x)$ respectively the closed
and the open balls with radius $\rho$ centered at $x$. A set $\mathcal{A}\subset M$ 
is a {\it Delone set}, if for some $0<r<R$ the following conditions hold.
\begin{itemize}

\item $B_r^\circ(x)\bigcap B_r^\circ(y)=\varnothing$ for every $x,y\in\mathcal{A}$.

\item $\bigcup\limits_{x\in \mathcal{A}} B_R(x)=M$.

\end{itemize}

Two Delone sets $\mathcal{A}\subset M_1$ and $\mathcal{B}\subset M_2$ are {\it biLipschitz equivalent},
if there exist a real $\lambda\geqslant 1$ and a bijection $F:\mathcal{A}\to\mathcal{B}$ such that the inequality 
$$\frac 1\lambda d_{M_1}(x,y)\leqslant d_{M_2}(F(x),F(y))\leqslant \lambda d_{M_1}(x,y)$$
holds for every $x,y\in \mathcal{A}$.
 
Map $F$ for which such an inequality holds is called $\lambda$-{\it biLipschitz}.

The question about biLipschitz equivalence was raised by M.~Gromov in [1]. In particular,
the following problem was stated:

{\it Given a metric space $M$ determine if every two Delone sets 
$\mathcal A, \mathcal B \subset M$ are biLipschitz equivalent.}

If $M=\mathbb E^1$ --- a Euclidean line, then the answer is obviously positive. Also positive answers were obtained by
P.~Papasoglu (see [2]) for homogeneous trees, O.~Bogopolsky (see [3]) for hyperbolic spaces $\mathbb H^d$
and K.~Whyte (see [4]) for non-amenable spaces.

In case of Euclidean space $M=\mathbb E^d$ of dimension $d\geqslant 2$ D.~Burago and B.~Kleiner (see [5])
and independently C.McMullen (see [6]) proved the following result:

{\bf Theorem 1.} {\it For every integer $d\geqslant 2$ there exists a Delone set $\mathcal{A}\subset \mathbb E^d$
which is not biLipschitz equivalent to the integer net $\mathbb{Z}^d$.}

In~[5] theorem 1 is proved for $d=2$, but the proof is easily generalized for every dimension
$d\geqslant 2$. Therefore in $\mathbb E^d$ for every $d\geqslant 2$ there exist at least 2 biLipschitz classes.

The main result of this paper is

{\bf Theorem 2.} {\it For every integer $d\geqslant 2$ the set of biLipschitz classes in $\mathbb{E}^d$ has cardinality continuum.}

\medskip

{\bf Proof of theorem 2}
 
\medskip

Obtain the upper estimate for cardinality of the set of biLipschitz classes.
 
Use the following result of A.~Garber (see [7, lemma 2]) 

{\bf Lemma 3.} {\it Let $\mathcal A\subset \mathbb E^d$ be a Delone set. Then there exists a Delone set
$\mathcal D\subset \mathbb Z^d$ such that $\mathcal A$ and $\mathcal D$ are biLipschitz equivlent.}

From lemma 3 follows that every biLipschitz class has at least one member among subsets of $\mathbb Z^d$.
Therefore cardinality of the set of biLipscitz classes is at most cardinality of family of all subsets containing in $\mathbb Z^d$, 
i.e. continuum. The upper estimate proved.

To prove the lower estimate obtain a continuum family of pairwise non-equivalent Delone sets. These sets will be members of some special class.

From this point we consider only rectangular coordinates in $\mathbb E^d$. Parallelepipeds (cubes) with edges parallel to
coordinate lines are called {\it coordinate}.

Let $Q$ be a coordinate cube. Denote by $m(Q)$ its vertex with the least sum of coordinates.

Consider a tiling $T$ of $\mathbb E^d$ into coordinate cubes whose edge lengths belong to $[1,L]$. 
The set $\mathcal A=\{m(Q):Q\in T\}$ is obviously
a Delone set. Delone sets obtained in such a way are called $L$-{\it special}.

Consider a map $G_{\mathcal A}:\mathcal A\to T$ sending each point $x\in \mathcal A$ to a cube $Q\in T$ such that $x=m(Q)$.

A point $x$ of special Delone set $\mathcal A$ is {\it standard} if $G_{\mathcal A}(x)$ a unit cube and 
{\it exceptional} otherwise.

Introduce some notation. Let 
$$\mathcal P_{MN}=\mathbb Z^d\cap ([0, MN]\times[0,N)^{d-1}),$$ 
$$\mathcal P^i_{MN}=\mathbb Z^d\cap ([iN, (i+1)N)\times[0,N)^{d-1}) \;\text{for}\; i=0,1,\ldots, M-1,$$ 
$$u=(0,0,\ldots, 0),$$ 
$$v=(MN, 0, 0, \ldots, 0).$$ 

Call points $x,y\in\mathcal P_{MN}$ {\it corresponding}
if $y-x=(N,0,0,\ldots,0)$.

{\bf Lemma 4.} {\it Let $\lambda\geqslant 1$, $\varepsilon\in (0,\frac 14)$ è $a\in (0,1)$. Then there exist 
$k>0$ and $M_0\in\mathbb N$ such that for every $M,N\in\mathbb N$ ñ $M>M_0$ and for arbitrary $\lambda$-biLipschitz map 
$F: \mathcal{P}_{MN} \to \mathbb{E}^d$ at least one of the following statements hold:
\begin {enumerate}

\item There exist corresponding points $x,y$ such that
$$\frac{|F(y)-F(x)|}{|y-x|}>(1+k)\frac{|F(v)-F(u)|}{|v-u|};$$

\item There exists $i$ such that number of pairs of corresponding points $x\in \mathcal P^i_{MN}$, $y\in\mathcal P^{i+1}_{MN}$
for which holds
$$\frac{|F(y)-F(x)-\frac 1M (F(v)-F(u))|}{\frac 1M |F(v)-F(u)|}<\varepsilon,$$
is at least $aN^d$.

\end {enumerate}}

Proof for $d=2$ is in [5, Lemma 3.2]. Proof for an arbitrary $d$ is obtained by a straightforward repeating the arguments of [5].

{\bf Lemma 5.} {\it Let $I=[0,1]$, $\alpha\in(0,\frac 12)$, and let $P,Q\subset I^d$ be closed sets 
with a boundary being a finite polyhedron. If $P\cup Q=I^d$, $\inter P\cap \inter Q=\varnothing$ and also
$\vol_d (P)\geqslant\alpha$ and $\vol_d (Q)\geqslant\alpha$ then $(d-1)$-dimensional volume 
$\vol_{d-1} (\partial P\cap \partial Q)\geqslant \frac {\alpha}{2^{d-1}}$.}

\textbf{Proof.} Denote by $\pi$ the projection onto hyperplane \newline $x_1=0$.

Conduct the proof by induction over $d$. 

Induction base: $d=2$. If $\vol_1 (\pi(P)\cap\pi(Q))\geqslant \frac\alpha 2$ then statement of lemma is obviously true. Otherwise 
the following inequalities hold: 
$$1-\alpha\geqslant 1-\vol_2(Q)=\vol_2 (P)\geqslant \vol_1(\pi(P))-\vol_1 (\pi(P)\cap\pi(Q)).$$ 
Hence $\vol_{d-1}(\pi(P))<1-\frac \alpha 2$.

Therefore there exists $t_P\in(0,1)$ such that $P\cap\{x_2=t_P\}=\varnothing$. Similarly, there exists $t_Q\in(0,1)$ such that
$Q\cap\{x_2=t_Q\}=\varnothing$.
It follows that projection of $P,Q$ onto line $x_2=0$ is a segment $[0,1]$. Hence 
$\vol_1 (\partial P\cap \partial Q)\geqslant 1>\frac \alpha 2$ and for $d=2$ statement is proved.

Induction step. Similarly to previous if $\vol_{d-1} (\pi(P)\cap\pi(Q))\geqslant \frac\alpha 2$ statement of lemma is obvious. Otherwise
$\vol_{d-1}(\pi(P))<1-\frac \alpha 2$. Then every section of $Q$ by a hyperplane $x_1=t$ has a $(d-1)$-volume $\geqslant\frac \alpha 2$.
Similarly, every section of $P$ by a hyperplane $x_1=t$ has a $(d-1)$-volume $\geqslant\frac \alpha 2$. By induction assumption, every
section of $\partial P\cap \partial Q$ has a  $(d-2)$-volume $\geqslant\frac {\alpha} {2^{d-1}}$, and the statement of lemma is now obvious.

{\bf Lemma 6.} {\it   Given $\lambda\geqslant 1$, $L\geqslant 1$ and rational $c>1$ there exists a finite point set 
$\mathcal B_0$ and a parallelepiped $\Pi=\prod\limits_{i=1}^d [0, b_i)$ ñ $b_i\in \mathbb N$ such that:
      
      \begin{enumerate}
      
         \item $\mathcal B_0\subset\Pi$.
         
         \item There exists a tiling $T_0$ of $\Pi$ into coordinate cubes with edges 1 and $c$ such that 
         $\{m(Q):Q\in T_0\}=\mathcal B_0$.
            
         \item For every Delone set $\mathcal B$ such that $\mathcal B\cap \Pi=\mathcal B_0$ and \newline
            $(0,0,\ldots,0,b_d)\in\mathcal B$ and for every $\lambda$-biLipschitz bijection
            $F:\mathcal B\to\mathcal A$ where $\mathcal A$ is $L$-special, the set $F(\mathcal B_0)$ has at least one
            exceptional point.
      \end  {enumerate}}

{\bf Proof.} Conduct the construction of $\mathcal B_0$ in 3 steps:

\begin{enumerate}

\item Choose $\varepsilon$ è $a$ which have the same meaning as in lemma 4; choose a parameter $H_0$.

\item Choose $N$ and $M$.

\item Choose $H$ fulfilling $H\geqslant H_0$ and construction of $\mathcal B_0$ itself.

\end{enumerate}

Describe the construction beginning from the last step. Let $\varepsilon, a, M, N, H_0$ be already chosen on previous steps.

Take a parallelepiped
$$\Phi_{0,(0,0,\ldots, 0)}=[0,1)^{d-1}\times [0,M).$$ Consider its tiling into unit cubes. 
Colour these cubes checkerboardwise into black an white, starting with black.

Take in parallelepiped $\Phi_{0,(0,0,\ldots, 0)}$ parallelepipeds 
$$\Phi_{1,\frac 1N\cdot(j_1,j_2,\ldots, j_{d-1}, 0)}=[0,\frac 1M)^{d-1}\times [0,M)+\frac 1N\cdot(j_1,j_2,\ldots, j_{d-1}, 0)$$
where $j_i=0,1,\ldots, N-1$. From this point colouring of $\Phi_{0,(0,0,\ldots, 0)}$ will change only inside parallelepipeds 
of type $\Phi_{1,z}$. Divide each of these parallelepipeds into
cubes with edge equal to $\frac 1M$ and colour them checkerboardwise starting from black. 

Continue the process. On $\nu$-th step in each parallelepiped of type 
$$\Phi_{\nu-1, z}=[0, \frac 1{M^{\nu-1}})^{d-1}\times [0,M)+z$$
take the parallelepipeds
\begin{align*}
&\Phi_{\nu, z+\frac 1{N\cdot M^{\nu-1}}\cdot(j_1,j_2,\ldots, j_{d-1}, 0)}=\\
&=[0, \frac 1{M^\nu})^{d-1}\times [0,M)+z+\frac 1{N\cdot M^{\nu-1}}\cdot(j_1,j_2,\ldots, j_{d-1}, 0).
\end{align*}
From this point colouring will change only inside these parallelepipeds. Divide each of these parallelepipeds into
cubes with edge equal to $\frac 1{M^l}$
and colour them checkerboardwise starting from black.

Repeat while $\nu\leqslant\nu_0= \lceil \log_{1+k}\lambda^2\rceil+2$.

Note that if $\Phi_{0,(0,0,\ldots, 0)}$ is divided into cubes with edge $\frac 1 {NM^{\nu_0}}$ then each
of them is coloured in one colour --- black or white. Call them {\it coloured cubes}

Make a homothety of parallelepiped $\Phi_{0,(0,0,\ldots, 0)}$ together with colouring of coefficient $H$ and center at origin. Choose
$H$ such that coloured cubes were taken into cubes that have integer edges and also could be divided into cubes with edge $c$.
Inequality $H\geqslant H_0$ also must hold.

Images of black coloured cubes divide into unit cubes and images of white cubes --- into cubes with edge $c$.
The obtained tiling of $\Pi=H\cdot \Phi_{0,(0,0,\ldots, 0)}$ denote by $T_0$. Let $\mathcal B_0=\{m(Q):Q\in T_0\}$.

Describe the second step. Let $\varepsilon$, $a$ è $H_0$ be chosen before, choose $N$ and $M$.
 
Let $y-x=(0,0,\ldots,0,1)$. Choose $N$ such that if $P=[0,\frac 1N]^d$ then for every vector $\mathbf{e}$
fulfilling
$$|F(y)-F(x)-\mathbf{e})|<\varepsilon \cdot|\mathbf{e}|,$$
holds the inequality
$$|F(y')-F(x')-\mathbf{e}|<2\varepsilon \cdot|\mathbf{e}|$$
if only $x'\in x+P$, $y'\in y+P$ and $F$ is $\lambda$-biLipschitz.

Let $P_1$ è $P_2$ be cubes with edge $\frac H{M^l}$ coloured on $l$-th step black and white respectively.
Let each be divided into $N^d$ equal cubes and let $Q_1\subset P_1$ è $Q_2\subset P_2$ be
such cubes. Choose $M$ such that independently from choice of $H$ holds true 
$$\frac{\#(P_1\cap\mathcal B_0)}{\#(P_2\cap\mathcal B_0)}\geqslant\frac{1+c}2.$$

This inequality is obviously true if only
$$\vol_d \left(Q_1 \cap (\cup_z \Phi_{l+1,z}) \right)\leqslant \frac 1{c-1} \vol_d(Q_1)\;\text{and}$$
$$\vol_d \left(Q_2 \cap (\cup_z \Phi_{l+1,z}) \right)\leqslant \frac 1{c-1} \vol_d(Q_2)$$
which is true for big enough $M$. Also take $M>M_0$ where $M_0$ comes from lemma 4 and $N|M$.

Describe the first step.

Let $a=\frac{3+c}{2+2c}$. Show that there exists a choice of $\varepsilon$ and $H_0(\varepsilon)$ such that
$\mathcal B_0$ constructed as before fulfilled the conditions of lemma 6. 

Suppose that for every $\varepsilon$ and $H_0$ there is a Delone set $\mathcal B\supset\mathcal B_0$ fulfilling the 
conditions of lemma 6 and $\lambda$-biLipschitz bijection
$F: \mathcal B \to \mathcal A$ such that $F(\mathcal B_0)$ consists only of standard points. 

Let $u=(0,0,\ldots,0)$, $v=(0,0,\ldots,0,HM)$. If the first case of statement of lemma 4 holds 
there exist corresponding points $x,y$ such that \newline
$\frac{|F(y)-F(x)|}{|y-x|}>(1+k)\frac{|F(v)-F(u)|}{|v-u|}$. In this case instead of $u,v$ consider a pair $x,y$ and  restriction of 
$F$ to a subset of $\mathcal B_0$ contained in parallelepiped 
$$x+ [0,\frac HM)^{d-1}\times [0,H).$$
Apply to this set all the arguments similarly as to $\mathcal B_0$. If such a substitution can be made $\lceil \log_{1+k}\lambda^2\rceil+2$ times,
then from $u,v$ we come to $u', v'$ such that 
$\frac{|F(v')-F(u')|}{|v'-u'|}>\lambda^2\frac{|F(v)-F(u)|}{|v-u|}$, which makes a contradiction to $\lambda$-biLipschitz property of $F$.

Therefore on some step we have the second case of lemma 4. Let the adjoint cubes for which this case holds
have numbers $i$ è $i+1$. Let also $i$-th cube be originally white and, respectively, $(i+1)$-th black.
Let $\tilde F=G_{\mathcal A}\circ F$.

Let $\mathcal{C}$ be a set of points of $i$-th cube such that have distance at least $10\lambda L$ from its boundary, $\mathcal C'$ 
are all points of $i+1$-th cube. Using our assumptions 
obtain two inequalities involving $\vol_{d-1}(\partial\tilde F(\mathcal{C}))$. 
$$\vol_{d-1}(\partial\tilde F(\mathcal{C}))\leqslant\beta_0\cdot |\mathcal{C}|^{\frac {d-1}d},$$
$$\vol_{d-1}(\partial\tilde F(\mathcal{C}))\geqslant\beta_2\varepsilon^{-1}\cdot |\mathcal{C}|^{\frac {d-1}d}.$$
For small enough $\varepsilon$ they contradict each other and that completes the proof of lemma 6.

Proving lemmas 7 and 8 $H$ is assumed big enough depending on $\varepsilon$,
i.e. $H\geqslant H_0(\varepsilon)$.

{\bf Lemma 7.} {\it Inequality
$$\vol_{d-1}(\partial\tilde F(\mathcal{C}))\leqslant\beta_0\cdot |\mathcal{C}|^{\frac {d-1}d}$$
holds true, where $\beta_0$ is a constant depending on $d$, $\lambda$, $L$ and $c$ (but not $\varepsilon$).}

{\bf Proof.} This inequality follows immediately from the fact that if $\tilde F(x)$ has common boundary with $\tilde F(\mathcal C)$
then $x$ is a point of $i$-th cube, 
that does not depend on $\mathcal C$. The number of such points does not exceed $\beta_1\cdot |\mathcal{C}|^{\frac {d-1}d}$, hence $(d-1)$-volume
of boundary of corresponding cubes does not exceed $\beta_0\cdot |\mathcal{C}|^{\frac {d-1}d}$.

Let $s=4\varepsilon\frac 1M |F(v)-F(u)|$.

Let $K$ be a real independent from $\varepsilon$ and such that 
$$\left(1+\frac{(2K+2)^d-(2K)^d}{K^d}\right)\frac 4{3+c}<\frac 8{7+c}.$$
 
Take a full (in respect to inclusion relationship) packing of coordinate cubes with centers in $F(\mathcal{C})$ and edges equal to $Ks$.
Let it consist of $W$ cubes.
 
Denote by $U$ a union of coordinate cubes with the same centers and edges equal to $2Ks$. 
Since the chosen packing is full all points of $F(\mathcal{C})$ are contained in $U$.

Let $\tau$ be a translation by vector $\frac 1M (F(v)-F(u))$.   
Consider $\frac s2$-neighbourhood of $\tau(F(\mathcal{C}))$. Denote by $\mathcal C''$ the set of points of $\mathcal A$ that belong to
this neighbourhood. Since the second case of lemma 4 assumed true and due to choice of $a$ and $M$ obtain: 
$$|\mathcal C'\cap \mathcal C''|\geqslant a\frac{1+c}2 \cdot|\mathcal{C}|=\frac{3+c}4 \cdot|\mathcal{C}|.$$ 
But $\tilde F(\mathcal C'\cap \mathcal C'')$ is contained in the union of cubes with the same centers as $\tau(U)$ 
and edge equal to $(2K+2)s$, because for big $H$ holds $s>1$. Denote this union by $U_1$.

Note that $\vol_d (U)\geqslant W K^d s^d$, $\vol_d(U_1)\leqslant \vol_d(U)+((2K+2)^d-(2K)^d)Ws^d$. According to choice of $K$ obtain
$$\vol_d (\tilde F(\mathcal C))\leqslant \frac 4{3+c} \vol_d(U_1) \leqslant \frac 8{7+c}\vol_d(U).$$ 

Rewrite the last inequality as
$\vol_d(U\setminus\tilde F(\mathcal C))\geqslant\frac {c-1}{7+c}\vol_d(U)$. Due to an estimate for $\vol_d (U)$ already obtained,
$$\vol_d(U\setminus\tilde F(\mathcal C))\geqslant\frac {c-1}{7+c}W K^d s^d.$$

Let $\mu\in(0,1)$ be such that $\mu+(1-\mu)\cdot\frac{c-1}{2^d(10+c)}<\frac{c-1}{2^d(7+c)}$. Note that $\mu$ does not depend on
$\varepsilon$. Then in at least $\mu W$ cubes of $U$ set $\tilde F(\mathcal C)$ occupies volume at most 
$\left(1-\frac{c-1}{2^d(10+c)}\right)\cdot (2Ks)^d$. Call the cubes {\it marked}.

{\bf Lemma 8.} {\it Suppose $\varepsilon$ small enough, then in our assumptions on $F$
$$\vol_{d-1}(\partial \tilde F(\mathcal{C}))\geqslant\beta_1\varepsilon^{-1}\cdot |\mathcal{C}|^{\frac {d-1}d},$$
where $\beta_2$ depends on $d$, $\lambda$, and $c$.}

{\bf Proof.}  If $\varepsilon<\frac1{4K\lambda^2}$ and $H$ is big enough then in every marked cube 
$\tilde F(\mathcal C)$ occupies volume at least $\beta_3\cdot(2Ks)^d$. Indeed, if $F(x)$ is a center of marked cube
then due to $\lambda$-biLipschitz property of $F$ all points of $\mathcal C\cap B_{\frac {Ks}\lambda}(x)$ are taken inside this cube.
For $H$ big enough $s$ is also big, then 
$$\vol_d (\tilde F(\mathcal C\cap B_{\frac {Ks}\lambda}(x)))=|\mathcal C\cap B_{\frac {Ks}\lambda}(x))|\geqslant 2\beta_3 \cdot(2Ks)^d,$$
and on the other hand, a part of volume of $\tilde F(\mathcal C\cap B_{\frac {Ks}\lambda}(x))$ not exceeding\newline
$(2Ks+2)^d-(2Ks)^d$ can be excluded from the marked cube. But for big enough $s$ it does not exceed $\beta_3 \cdot(2Ks)^d$. Hence the inequality.

According to lemma 5 inside marked cubes 
$\partial \tilde F(\mathcal C)$ has $(d-1)$-volume at least $\beta_4 s^{d-1}$ where $\beta_4$ depends on $d$, $\lambda$, and $c$.

Since cubes of packing do not intersect, no $8^d+1$ cubes of $U$ have a common point. 
Then for some $\beta_5$, depending on $d$, $\lambda$, and $c$ holds 
$$\vol_{d-1}(\partial\tilde F(\mathcal C))\geqslant \beta_5 Ws^{d-1}.$$

Since $F(\mathcal C)\subset U$ obtain $|\mathcal C|\leqslant W(2Ks+2)^d\leqslant \beta_6 W s^d$. Again $H$ and $s$ are assumed big enough. Therefore
$$\vol_{d-1} (\partial\tilde F(\mathcal C))\geqslant \beta_7 |\mathcal C| s^{-1}.$$

Due to $\lambda$-biLipschitz property of $F$ holds 
$$\frac 1M |F(v)-F(u)|\leqslant \beta_8 |\mathcal C|^{\frac 1d}.$$ 

Using the definition of $s$ obtain
$$s\leqslant\beta_9\varepsilon|\mathcal C|^{\frac 1d},$$
which together with the last inequality for $\vol_{d-1} (\partial\tilde F(\mathcal C))$ implies the statement of lemma.

{\bf Lemma 9.} {\it   Given real $\lambda\geqslant 1$, $L\geqslant 1$, rational $c>1$ and positive integer $j\in\mathbb N$ there exists a finite
point set $\tilde {\mathcal D}$ and a parallelepiped $\Pi=\prod\limits_{i=1}^d [0, b_i)$ where $b_i$ are positive integer such that:
      
      \begin{enumerate}
      
         \item $\tilde {\mathcal D}\subset\Pi$.
         
         \item There exists a tiling $T_0$ of $\Pi$ into coordinate cubes with edges 1 and $c$ such that
         $\{m(Q):Q\in T_0\}=\tilde {\mathcal D}$.
            
         \item For every Delone set $\mathcal D$ fulfilling $\mathcal D\cap \Pi=\tilde {\mathcal D}$ and \newline
            $(0,0,\ldots,0,b_d)\in\mathcal D$ for every $\lambda$-biLipschitz bijection
            $F:\mathcal D\to\mathcal A$ with $L$-special Delone set $\mathcal A$, the set $F(\tilde {\mathcal D})$ contains
            at least $j$ exceptional points.
      \end  {enumerate}}

{\bf Proof.} If $j=1$ then the desired statement is exactly lemma 6. If $j>1$ take a parallelepiped $\Pi(\lambda, L, c, j)$ with first $d-1$ 
edges equal to corresponding edges of  $\Pi(\lambda, L, c, 1)$ and the last edge $j$ times greater than the corresponding edge of $\Pi(\lambda, L, c, 1)$.
Divide $\Pi(\lambda, L, c, j)$ into $j$ parallelepipeds congruent to $\Pi(\lambda, L, c, 1)$. Take in each of them a set congruent to 
$\tilde {\mathcal D}(\lambda, L, c, 1)$. Denote the obtained set by $\tilde {\mathcal D}(\lambda, L, c, j)$. Obviously, it fulfills the statement of lemma.

Return to the proof of theorem 2. Let $\{c_i\}_{i=1}^\infty$ be a sequence of rationals from $(1,2]$, e.g.  $c_i=1+\frac 1i$.

In noataion of lemma 9 let $\mathcal{D}_1=\tilde {\mathcal D}(1,2,c_1,1)$. By induction define 
$$\mathcal{D}_j=\tilde{\mathcal D}(j,2,c_j,\sum\limits_{i=1}^{j-1}\#\mathcal{D}_i+1).$$ 
Let $r_j=100j\cdot\diam(\mathcal D_{j+1})$. Without loss of generality, let $r_j$ be strictly increasing.

Let $G$ be an additive group of rationals with denominator equal to some positive integer exponent of 2. From each class of
$\mathbb{R}/G$ choose one number and for each chosen number take the sequence of digits after the point in its binary representation. 
Obtain a continuum set of non-confinal $(0,1)$-sequences, i.e every two sequences have an infinite set of indices for which the corresponding members
are different.

For every taken sequence $\alpha=\{\alpha_i\}_{i=1}^\infty$ construct a 2-special Delone set $\mathcal{D}_\alpha$ as follows. Take $\mathcal{D}_1$ so that 
the corresponding parallelepiped was coordinate with integer vertices. Further, if $\alpha$ has zero as $j$-th digit take a copy of $\mathcal{D}_{j+1}$
at $\lceil r_j\rceil$ to the right from $\mathcal{D}_{j}$; if $\alpha$ has unit as $j$-th digit then take a copy of $\mathcal{D}_{j+1}$
at $\lceil 100jr_j\rceil$ to the right from $\mathcal{D}_{j}$. Also corresponding to $\mathcal{D}_{j+1}$ parallelepiped should be coordinate 
with integer vertices. In addition, include into $\mathcal{D}_\alpha$ all points of $\mathbb{Z}^d$ which are outside all the parallelepipeds
corresponding to $\mathcal{D}_j$. These points will be standard for $\mathcal{D}_\alpha$.

Prove that any two constructed sets are not biLipschitz equivalent.

Let there exist $\lambda$-biLipschitz bijection $F:\mathcal{D}_\alpha \to \mathcal{D}_\beta$. Lemma 10
states some property of this bijection.

{\bf Lemma 10.} {\it For $j>2\lambda$ there exists a point of $\mathcal{D}_\alpha$ in a copy of $\mathcal{D}_j$ (see construction of $\mathcal{D}_\alpha$)
which is sent by $F$ into some point of copy of $\mathcal{D}_j$ in $\mathcal{D}_\beta$.}

{\bf Proof.} Indeed, since image of $\mathcal{D}_j$ contains many enough exceptional points, there is a point $x$ of $\mathcal{D}_\alpha$ in a 
copy of $\mathcal{D}_j$ that is sent into an exceptional point, and, moreover, $F(x)$ does not belong to copies of $\mathcal{D}_i$ for $i<j$. If $F(x)$ 
belongs to copy of $\mathcal{D}_j$ then the proof is complete. Let $y=F(x)$ belong to copy of $\mathcal{D}_{j+k}$, $k>0$. Consider the image of
$\mathcal{D}_{j+k}$ under $F^{-1}$ that is also $\lambda$-biLipschitz. Let $z\in\mathcal{D}_{j+k}$. Then $$|F^{-1}(y)-F^{-1}(z)|\leqslant\lambda|y-z|\leqslant\lambda\diam(\mathcal{D}_{j+k})<r_{j+k-1}.$$ 
Therefore all exceptional points in image of $\mathcal{D}_{j+k}$ under $F^{-1}$ belong to copies of $\mathcal{D}_i$ for $i<j+k$, and it immediately 
implies a contradiction since there are many enough exceptional points in the image of $\mathcal{D}_{j+k}$.

Continue the proof of theorem 2.

Take $j>2\lambda$ such that $\alpha$ and $\beta$ differ in $j$-th digit. Without loss of generality, 
$\alpha_j=0$ and $\beta_j=1$. By lemma 10 there are $x_j$ and $x_{j+1}$ in copies of $\mathcal{D}_j$ and $\mathcal{D}_{j+1}$ respectively in 
$\mathcal{D}_\alpha$ which are taken to points $y_j$ and $y_{j+1}$ of corresponding copies in $\mathcal{D}_\beta$. By construction of $\mathcal{D}_\alpha$
and $\mathcal{D}_\beta$,
$$100j\cdot\diam(\mathcal{D}_{j+1})<|x_j-x_{j+1}|<(100j+2)\cdot\diam(\mathcal{D}_{j+1}),\; \text{and}$$  
$$10000j^2\cdot\diam(\mathcal{D}_{j+1})<|y_j-y_{j+1}|<(10000j^2+2)\cdot\diam(\mathcal{D}_{j+1}).$$
Therefore 
$$\frac{|y_j-y_{j+1}|}{|x_j-x_{j+1}|}>99j>\lambda.$$ 
A contradiction with $\lambda$-biLipschitz property of $F$ makes proof of theorem~2 complete.

\smallskip 
\centerline {\bf References} 

1.~M.~Gromov. Asymptotic invariants for infinite groups~// London
Mathematical Society Lecture Notes, vol.~182, Geometric group theory. eds.
J.~A.~Niblo, M.~A.~Roller, J.~W.~S.~Cassels, 1993.

2.~P.~Papasoglu. Homogeneous trees are bi-Lipschitz equivalent~// Geom.
Dedicata, vol.~54, 1995, 301-306.

3.~O.~V.~Bogopolskii, Infinite commensurable hyperbolic groups
are biLipschitz equivalent~// Algebra and Logic, vol.~36, no.~3, 1997,
155-163.

4.~K.~Whyte, Amenability, bi-Lipschitz equivalence, and the von Neumann
conjecture~// Duke Math. J. vol. 99, 1999, 93-112.

5.~D.~Burago, B.~Kleiner, Separated nets in Euclidean space and Jacobians of
bi-Lipschitz maps~//Geom. Funct. Anal. vol.~8, 1998, 273-282.

6.~C.~McMullen, Lipschitz maps and nets in Euclidean space~// Geom. Funct.
Anal. vol.~8, 1998, 304-314.

7.~A.~I.~Garber, On equivalence classes of separated nets~// Model.
and Anal. of Inf. Syst., vol.~16, no.~2, 2009, 109-118 (in Russian). 

\end  {document}